\def\draftdate{May 23, 2014}

\documentclass{amsart}
\usepackage{amssymb}
\usepackage{xy}

\newcommand{\Wh}[1]{Wh(#1)}

\newcommand{\trc}[1]{\mathit{trc}_{#1}}

\newcommand{\cmp}{\mathit{cmp}}
\newcommand{\cmppi}{\mathit{cmp}'}

\newcommand{\et}{\text{\'et}}

\def\myop#1{\mathop{\textstyle #1}\limits}
\let\iso\cong
\let\sma\wedge

\renewcommand{\to}{\mathchoice{\longrightarrow}{\rightarrow}{\rightarrow}{\rightarrow}}

\newcommand{\sto}{\rightarrow}

\newcommand{\phat}{^{\scriptscriptstyle\wedge}_{p}}

\let\catsymbfont\mathcal
\newcommand{\aD}{{\catsymbfont{D}}}
\newcommand{\aK}{{\catsymbfont{K}}}
\newcommand{\aM}{{\catsymbfont{M}}}
\newcommand{\aX}{{\catsymbfont{X}}}

\newcommand{\bA}{{\mathbb{A}}}
\newcommand{\bC}{{\mathbb{C}}}
\newcommand{\bCP}{\bC P}

\newcommand{\bS}{{\mathbb{S}}}
\newcommand{\bZ}{{\mathbb{Z}}}
\newcommand{\bQ}{{\mathbb{Q}}}

\newcommand{\SMot}{\aM ot_{\mathrm{ex}}}
\newcommand{\DMot}{\aD(\SMot)}

\def\quickop#1{\expandafter\DeclareMathOperator\csname
#1\endcsname{#1}}
\quickop{id}\quickop{Id}
\quickop{Tor}
\quickop{holim}\quickop{hocolim}\quickop{hofib}
\quickop{Spec}
\quickop{Gal}
\quickop{Hom}

\numberwithin{equation}{section}

\newtheorem{thm}[equation]{Theorem}
\newtheorem{main}{Theorem}

\newtheorem{lem}{Lemma}
\newtheorem{prop}[equation]{Proposition}
\theoremstyle{definition}

\theoremstyle{remark}

\xyoption{arrow}
\xyoption{matrix}
\xyoption{cmtip}
\SelectTips{cm}{}

\newcommand{\term}[1]{\textit{#1}}

\newdir{ >}{{}*!/-5pt/\dir{>}}

\begin{document}

\title[The nilpotence theorem for $K(\bS)$]
{The nilpotence theorem for the algebraic $K$-theory of the
sphere spectrum}

\author{Andrew J. Blumberg}
\address{Department of Mathematics, The University of Texas,
Austin, TX \ 78712}
\email{blumberg@math.utexas.edu}
\thanks{The first author was supported in part by NSF grant
DMS-1151577}
\author{Michael A. Mandell}
\address{Department of Mathematics, Indiana University,
Bloomington, IN \ 47405}
\email{mmandell@indiana.edu}
\thanks{The second author was supported in part by NSF grant
DMS-1105255}

\date{\draftdate} 

\begin{abstract}
We prove that in the graded commutative ring $K_{*}(\bS)$, all positive
degree elements are multiplicatively nilpotent.  The analogous
statements also hold for $TC_{*}(S;\bZ\phat)$ and $K_{*}(\bZ)$.  
\end{abstract}

\maketitle

\section{Introduction}

Much of the most exciting work in algebraic $K$-theory over the past
15 years has been aimed at the verification of the Quillen-Lichtenbaum
conjecture. The successful affirmation of this conjecture has led to
the identification of the homotopy types of the $K$-theory of the
integers $\bZ$ and the $K$-theory of the sphere spectrum $\bS$ at
regular primes~\cite{DwyerMitchell,RognesWeibel, Rognes2,Rognesp}. 
Since $H\bZ$ and $\bS$ are $E_{\infty}$ ring spectra,
$K(\bZ)$ and $K(\bS)$ are $E_{\infty}$ ring spectra and the graded
rings $K_*(\bS) = \pi_* K(\bS)$ and $K_*(\bZ) = \pi_* K(\bZ)$ are
commutative.  However, almost
nothing is known about the multiplicative structure.  The only work in
this direction so far is the investigation of Bergsaker and
Rognes~\cite{BergsakerRognes} of the Dyer-Lashof operations on
$TC_{*}(\bS)$ at the prime $2$.  The main theorem of this paper proves
the analogue of Nishida's Nilpotence Theorem for $K(\bS)$.

\begin{main}\label{thm:main}
Positive degree elements of $K_{*}(\bS)$ are nilpotent.   
\end{main}

On the way to proving the preceding theorem, we show the corresponding
nilpotence result for $K_*(\bZ)$.  We deduce this by observing that
$K_{2n(p-1)}(\bZ)\otimes \bZ_{(p)}=0$ for odd primes $p$ and $n>0$; it
can also be deduced from the multiplicative properties of the
Quillen-Lichtenbaum spectral sequence.

\begin{main}\label{thm:nilKZ}
Positive degree elements of $K_{*}(\bZ)$ are nilpotent.   
\end{main}

Much of the interest in $K(\bS)$ comes from its identification as
$A(*)$, Waldhausen's algebraic $K$-theory of the one-point space.
Work of Waldhausen and collaborators shows that $A(X)$ controls high
dimensional manifold theory (e.g., see~\cite{JahrenRognesWaldhausen}
and~\cite{WeissWilliamsOeuvre}) via the connection to the stable
pseudo-isotopy spectrum $\Wh{X}$.  Rognes shows that the infinite loop
space structure on $\Wh{*}$ that is relevant to the Hatcher-Waldhausen map
$G/O \to \Omega \Wh{*}$, where $G/O$ denotes the classifying
spectrum for smooth normal invariants, is induced by the ring
structure on $A(*)$~\cite{Rognes-HW}.  Moreover, $A(X)$ is a module
over $A(*)$; more generally, for any ring spectrum (or even any
Waldhausen category that admits factorization,
cf.~\cite{BlumbergMandellUW,BlumbergMandellTW}), the algebraic
$K$-theory spectrum is a module over $A(*)$.

Theorem~\ref{thm:main} also has direct implications in the context of
Kontsevich's non-commutative motives. The work of Blumberg, Gepner, and
Tabuada~\cite{BGT,BGT2} produces a candidate category of spectral
motives $\SMot$, which is a symmetric monoidal category with
objects the smooth and proper small stable idempotent-complete
$\infty$-categories.  The category of spectral motives is stable,
which in particular means that it has a tensor-triangulated homotopy
category and is enriched over spectra; the mapping spectra are
essentially bivariant algebraic $K$-theory.  The endomorphism spectrum
of the unit is precisely $K(\bS)$ (as an $E_{\infty}$ ring spectrum).  

The Devinatz-Hopkins-Smith nilpotence theorem and the Hopkins-Smith
thick subcategory theorem teach us that to understand a triangulated
category, we should look to its thick subcategories, which play the
role of prime ideals in derived algebraic
geometry~\cite{HopkinsGMH,NeemanChromatic,ThomasonCTS}.  More recently,
Balmer~\cite{BalmerTTC,BalmerSSS} proposes a systematic study of this
in the setting of ``tensor-triangulated geometry'', defining the
\term{triangulated spectrum} to be the space of prime proper thick
triangulated tensor ideals (with the Zariski topology).  Balmer
observes that there is a canonical map from the triangulated spectrum
to the spectrum of the graded ring of endomorphisms of the unit and
that in many known examples, the spectrum of the endomorphism ring
controls the triangulated spectrum of the tensor-triangulated
category.

Balmer~\cite{BalmerSSS} studies the example of the triangulated spectrum of (the
compact objects in) the Morel-Voevodsky stable $\bA^{1}$-homotopy
category $\aK^{\bA^{1}}$ (for a field $F$).  He uses the natural
composite from finite spectra to $\aK^{\bA^{1}}$ to finite spectra
(with the latter map given by the realization functor for $F\subset
\bC$) and the fact the algebraic Hopf map is non-nilpotent in the
endomorphisms of the unit in $\mathcal{K}^{\bA^{1}}$ to show that the
spectrum of $\aK^{{\bA^{1}}}$ is more complicated than the spectrum of
the stable category.  The derived category $\DMot$ of (dualizable
objects in) $\SMot$ admits a realization functor to finite spectra
induced by the Dennis trace map $K\to THH$ \cite[5.19]{BGT2}, which
splits the canonical map from finite spectra into $\DMot$ (since
$THH(\bS)\simeq \bS$).  Tabuada asked if the same argument could be
carried out for $\DMot$.  Our Theorem~\ref{thm:main} suggests that this kind of
approach cannot succeed when studying the triangulated spectrum of
non-commutative spectral motives.  Hesselholt has suggested that a
version of spectral motives based on ``real'' algebraic $K$-theory
(which has a $\bZ/2$-equivariant structure) might be more appropriate.

In a different direction, Morava has developed a conjectural program
for studying a homotopy-theoretic analogue of Kontsevich's
Grothendieck-Teichmuller group~\cite{Morava,KitchlooMorava} in terms
of homotopical descent for the category of spectral motives.  These
ideas revolve around understanding the structure of $\bS \sma^{\mathbf
L}_{K(\bS)} \bS$, which of course depends on the ring structure of
$K(\bS)$.  Morava notes that the calculation of this object is
straightforward rationally and results in a concise description as a
polynomial algebra on even degree generators: It is the polynomial
algebra on the free Lie co-algebra $L\langle
x_{6},x_{10},x_{14},\dotsc \rangle$ on generators in degrees 6, 10,
14, etc. (It is a Hopf algebra with coalgebra the tensor coalgebra on
$\pi_{*}\Sigma \Wh{*}_{\bQ}=\pi_{*}\Sigma^{6}ko_{\bQ}$, where $\Wh{*}$ is
the fiber of the map $K(S)\to S$.) Our results suggest that the
torsion part of this theory will be maximally complicated. 

\subsection*{Acknowledgments} The authors thank Bill Dwyer, Mike Hopkins, Michael Larsen, Tyler
Lawson, Barry Mazur, Jack Morava, and Justin Noel for helpful
conversations, and the IMA and MSRI for their hospitality while some
of this work was done. The authors especially thank Lars Hesselholt;
we learned much of what we know in this area through conversations
with him and from his unpublished notes~\cite{Lars-L} on $p$-adic
$L$-functions and the cyclotomic trace.

\section{Reduction of Theorems~\ref{thm:main} and~\ref{thm:nilKZ}}

Consider the arithmetic square
\[
\xymatrix@-1pc{%
K(\bS)\ar[r]\ar[d]&\myop\prod_{p}K(\bS)\phat\ar[d]\\
K(\bS)_{\bQ}\ar[r]&(\myop\prod_{p}K(\bS)\phat)_{\bQ}
}
\]
where $(-)\phat$ denotes $p$-completion (localization with respect to
the mod $p$ Moore spectrum) and $(-)_{\bQ}$ denotes rationalization.
To prove the nilpotence Theorem~\ref{thm:main}, it suffices to
prove the analogous nilpotence results for $K(\bS)_{\bQ}$ and
$K(\bS)\phat$ for each prime $p$; this is easy to see for $K(\bS)$
since $\pi_{*}K(\bS)$ is finitely generated in each degree
\cite[1.2]{DwyerAX}.  (Similar observations apply to $K(\bZ)$ for
Theorem~\ref{thm:nilKZ}.)  The rational part is well understood: The
natural map $K(\bS)_{\bQ} \to K(\bZ)_{\bQ}$ is an equivalence
\cite[2.3.8]{WaldhausenKT}, and classical results of
Borel~\cite[12.2]{Borel-1974} imply that the positive degree elements of
$\pi_* K(\bZ)_{\bQ}$ are concentrated in odd degrees and therefore
square to zero.  It remains to study the situation after $p$-completion.

Our strategy for studying the multiplicative structure on $K(\bS)\phat$
uses the cyclotomic trace map, which is a map of $E_\infty$ ring
spectra from $K(\bS)$ to the topological cyclic homology $TC(\bS)$.
The homotopy type of $TC(\bS)\phat$ (as a spectrum) is known by work
of~\cite{BHM}.

\begin{thm}[{\cite[5.16]{BHM}}]\label{thm:CP}
There is an equivalence of $p$-complete spectra
\[
TC(\bS)\phat\simeq \bS\phat \vee 
\hofib(\Sigma (\Sigma^{\infty}_{+}\bCP^{\infty}) \rightarrow \bS)\phat 
\simeq \bS\phat \vee (\bCP^{\infty}_{-1})\phat.
\]
\end{thm}

The Devinatz-Hopkins-Smith nilpotence theorems provide a criterion for
determining when elements in the homotopy groups of a ring spectrum
$R$ are multiplicatively nilpotent.  Specifically, an element $x\in
\pi_{*}R$ is nilpotent if and only if the Hurewicz map takes it to a
nilpotent element of $K(n)_{*}R$ for all $0\leq n\leq \infty$ (and all
primes $p$).  Although the previous theorem only identifies the
homotopy type of the underlying spectrum and says nothing about the
multiplication, it is enough to deduce a nilpotence result for
$TC(\bS)\phat$.

\begin{prop}\label{prop:tcodd}
Let $p$ be a prime, $0\leq n\leq \infty$, and let
$\widetilde{TC}(\bS;p)$ be the homotopy fiber of the augmentation
map $TC(\bS)\phat \to \bS\phat$ (obtained from the canonical map
$TC(\bS)\phat\to THH(\bS)\phat\simeq \bS\phat$).  Then
$K(n)_{*}(\widetilde{TC}(\bS;p))$ is concentrated in odd degrees.
\end{prop}

\begin{proof}
As a consequence of Theorem~\ref{thm:CP}, $\widetilde{TC}(\bS;p)\simeq
\Sigma(\bCP^{\infty}_{-1})\phat$.  The spectrum $\bCP^{\infty}_{-1}$
is the Thom spectrum of the virtual bundle $-\gamma$, for $\gamma$ the
tautological line bundle over $\bCP^{\infty}$.  The spectra $K(n)$ are
all complex oriented; the proposition now follows from
the Thom isomorphism.
\end{proof}

Since $\pi_{*}(TC(\bS)\phat)$ splits as $\pi_{*}\bS\phat\oplus
\pi_{*}\widetilde{TC}(\bS;p)$ with the first factor the image of
the inclusion of the unit, we obtain the following as an immediate
corollary of the previous proposition and the nilpotence theorem.

\begin{thm}\label{thm:tcodd}
For any prime $p$, all the non-zero degree elements of $\pi_{*}TC(\bS)\phat$
are nilpotent.
\end{thm}

In light of the previous result, Theorem~\ref{thm:main} becomes an
immediate consequence of the following lemma.  We prove this lemma for
odd $p$ in later sections; for $p=2$ it is a special case
of~\cite[3.16]{Rognes2}. 

\begin{lem}\label{lem:main}
For $p=2$, let $d=8$, and for $p$ odd, let $d=2(p-1)$.
The homotopy fiber of the cyclotomic trace map $\trc{p}\colon
K(\bS)\phat \to TC(\bS)\phat$ has trivial homotopy groups in degrees
$kd$ for $k>0$.
\end{lem}

\begin{proof}[Proof of Theorem~\ref{thm:main} from Lemma~\ref{lem:main}]
Given $x\in \pi_{k}(\bS\phat)$, $x^{d}\in
\pi_{kd}K(\bS\phat)$.  When $k>0$, we then know that for some power
$n$, $(x^{d})^{n}$ maps to zero in $\pi_{kdn}(TC(\bS)\phat)$
under the trace map by Theorem~\ref{thm:tcodd}.  By
Lemma~\ref{lem:main}, the kernel of the trace is zero in degree
$kdn$, and so $x^{kdn}=0$.
\end{proof}

As we used in the proof, Lemma~\ref{lem:main} implies that the
cyclotomic trace $K(\bS) \to TC(\bS)$ is injective in certain degrees.
In fact, for odd regular primes, the cyclotomic trace is injective in
all degrees.  This follows fairly easily from the work of Rognes on
$\Wh{*}$ at odd regular primes, specifically~\cite[3.6,3.8]{Rognesp}.
In the case of irregular primes, we expect that the trace fails to be
injective; we hope to return to this question in a future paper.

On the way to proving Lemma~\ref{lem:main}, we also prove the
following lemma.  It is well known that $\pi_{4k}(\bZ)\otimes
\bZ_{(p)}=0$ at regular primes, including $p=2$ (see, for example,
\cite[10.1]{Weibel-KBook}), and this combined with the following lemma
now proves Theorem~\ref{thm:nilKZ}.

\begin{lem}\label{lem:nilKZ}
For $p$ an odd prime, $\pi_{2(p-1)k}K(\bZ)\otimes \bZ_{(p)}=0$ for
$k>0$. 
\end{lem}

\section{Reduction of Lemmas~\ref{lem:main} and~\ref{lem:nilKZ}}

The basic strategy for the proof of Lemmas~\ref{lem:main} and
\ref{lem:nilKZ} is to
reduce the study of the homotopy fiber of the cyclotomic trace
$K(\bS) \to TC(\bS;p)$ to the study of the $p$-completion map
$\bZ[1/p] \to \bQ\phat$ in \'etale cohomology.  (This is now a fairly
standard approach; for instance, see~\cite[\S2-3]{Rognesp},
\cite{Geisser-trace}, \cite{Lars-L}.)  As
indicated above, from here on we assume that $p$ is odd (though all of
what we say would also apply in the case $p=2$ until~\eqref{eq:tss}).
First, we apply Dundas' theorem~\cite{Dundas} about the cyclotomic
trace: the square
\[
\xymatrix{
K(\bS)\phat \ar[r] \ar[d]_{\trc{p}} & TC(\bS; p) \ar[d]^{\trc{p}^{\bZ}} \\
K(\bZ)\phat \ar[r] & TC(\bZ; p)
}
\]
is homotopy cocartesian, where the horizontal maps arise from
linearization.  As a consequence, we have the following lemma:

\begin{prop}[Dundas~{\cite{Dundas}}]
The induced map $\hofib(\trc{p}) \to \hofib (\trc{p}^{\bZ})$ is an
equivalence.
\end{prop}

To understand $\hofib(\trc{p}^{\bZ})$, consider the commutative diagram
\[
\xymatrix@C+2ex{%
K(\bZ)\phat\ar[r]^{\cmp}\ar[d]_{\trc{p}^{\bZ}}
&K(\bZ\phat)\phat\ar[d]^{\trc{p}^{\bZ\phat}}\\
TC(\bZ)\phat\ar[r]_{\cmp^{TC}}&TC(\bZ\phat)\phat
}
\]
where the horizontal maps $\cmp$ and $\cmp^{TC}$ are induced by the
map of rings $\bZ\to \bZ\phat$. By work of
Hesselholt-Madsen~\cite{HM2}, the bottom map is a weak
equivalence~\cite[Addendum~5.2]{HM2} and the right hand map induces a weak
equivalence~\cite[Theorem~D]{HM2}  
\begin{equation}\label{eq:trcbzph}
K(\bZ\phat)\phat\to TC(\bZ\phat)\phat[0,\infty)
\end{equation}
(where $[0,\infty)$ denotes the
connective cover).  Thus, up to passing to a connective cover, we can
identify the trace map $\trc{p}^{\bZ}$ as the map $\cmp\colon
K(\bZ)\phat\to K(\bZ\phat)\phat$.  We then have the following
relationship between $\hofib(\trc{p})\simeq \hofib(\trc{p}^{\bZ})$ and
$\hofib(\cmp)$.

\begin{prop}
There is a cofiber sequence 
\[
\hofib(\cmp)\to \hofib(\trc{p})\to \Sigma^{-2}H\bZ\phat\to \Sigma \dotsb 
\]
\end{prop}

\begin{proof}
Using the equivalence of $\hofib(\trc{p})$ and $\hofib(\trc{p}^{\bZ})$
above, we get a diagram of cofiber sequences
\[
\xymatrix{%
\hofib(\cmp)\ar[r]\ar[d]
&K(\bZ)\phat\ar[r]^{\cmp}\ar@{=}[d]
&K(\bZ\phat)\phat\ar[r]\ar[d]_{\trc{p}^{\bZ\phat}}
&\Sigma \hofib(\cmp)\ar[d]\\
\hofib(\trc{p})\ar[r]
&K(\bZ)\phat\ar[r]&TC(\bZ\phat)\phat\ar[r]
&\Sigma \hofib(\trc{p})
}
\]
identifying, the right hand square as homotopy (co)cartesian.
Since $\pi_{-1}TC(\bZ)\phat=\bZ\phat$
and $\pi_{n}TC(\bZ)\phat=0$ for $n<-1$, the homotopy cofiber of the map
$\trc{p}^{\bZ\phat}$ in the diagram is $\Sigma^{-1}H\bZ\phat$.
Desuspending, we see that the homotopy cofiber of 
\[
\hofib(\cmp)\to \hofib(\trc{p})
\]
is $\Sigma^{-2}H\bZ\phat$.
\end{proof}

For Lemma~\ref{lem:main} then $\hofib(\cmp)$ works just as well as
$\hofib(\trc{p})$.  Quillen's localization sequence~\cite{Quillen-I}
gives cofiber sequences 
\begin{equation}\label{eq:QuillenLoc}
\begin{gathered}
\xymatrix{%
K(\bZ/p)\ar[r]\ar[d]_{\id}&K(\bZ)\ar[r]\ar[d]&K(\bZ[1/p])\ar[r]\ar[d]&\Sigma\dotsb \\
K(\bZ/p)\ar[r]&K(\bZ\phat)\ar[r]&K(\bQ\phat)\ar[r]&\Sigma\dotsb
}
\end{gathered}
\end{equation}
from which we can see that $\hofib(\cmp)$ is equivalent to the
homotopy fiber of the map
\[
\cmppi \colon K(\bZ[1/p])\phat\to K(\bQ\phat).
\]

\begin{prop}
There is a homotopy equivalence $\hofib(\cmp)\to\hofib(\cmppi)$.
\end{prop}

The advantage of this approach is that \'etale cohomology
methods at the prime $p$ can be applied in rings where $p$ is a unit.
Letting $R$ denote either $\bZ[1/p]$ or $\bQ\phat$, then $R$ satisfies
the ``mild extra hypotheses'' of 
Thomason~\cite[0.1]{ThomasonEtale}, which 
gives a spectral sequence
\begin{equation}\label{eq:tss}
E^{s,t}_{2}=H^{s}_{\et}(\Spec
R;\bZ/p^{n}(t/2))\quad\Longrightarrow\quad
\pi_{t-s}(K_{\et}(R);\bZ/p^{n})
\end{equation}
from \'etale cohomology to the mod $p^{n}$ homotopy groups of
(Dwyer-Friedlander) \'etale $K$-theory.  In the formula above
\[
\bZ/p^{n}(t/2)=\begin{cases}
\mu_{p^{n}}^{\otimes (t/2)}&\text{if }t\text{ is even}\\
0&\text{if }t\text{ is odd}
\end{cases}
\]
where $\mu_{p^{n}}$ denotes the $p^{n}$-th roots of $1$
(i.e., $\mu_{p^{n}}(A)=\{x\in A\mid x^{p^{n}}=1\}$, a sheaf in
the \'etale topology).  In this case the affirmed
Quillen-Lichtenbaum conjecture~\cite{Citation_needed} identifies
\[
\pi_{*}(K(R);\bZ/p^{n})=\pi_{*}(K_{\et}(R);\bZ/p^{n})
\]
for $*\geq 2$.  Also, because we have assumed that $p$ is odd,
$H_{\et}^{*}(R;\bZ/p^{n}(k))=0$ for $*>2$ \cite[\S
III.1.3]{Soule-KtheoryNumRing}, and the spectral sequence collapses to
give an isomorphism and a short exact sequence
\begin{equation}\label{eq:finKvsH}
\begin{gathered}
\pi_{2k-1}(K(R);\bZ/p^{n})\iso H^{1}_{\et}(\Spec R;\bZ/p^{n}(k))\\
0\sto H^{2}_{\et}(\Spec R;\bZ/p^{n}(k+1))\sto \pi_{2k}(K(R);\bZ/p^{n})
\sto H^{0}_{\et}(\Spec R;\bZ/p^{n}(k))\sto 0
\end{gathered}
\end{equation}
for $k>1$.  In fact, the calculation of the $H^{0}_{\et}$ term is well
known: 

\begin{prop}\label{prop:BernoulliDenom}
Let $R=\bZ[1/p]$ or $\bQ\phat$.  Then $H^{0}_{\et}(\Spec
R;\bZ/p^{n}(k))=0$ unless $(p-1)\mid k$ and for $k=m(p-1)$, $H^{0}_{\et}(\Spec
R;\bZ/p^{n}(k))\iso \mu^{\otimes k}_{p^{i}}(\bar Q)$, where
$p^{i}=\gcd(|m|p,p^{n})$ (for $m\neq 0$ or $i=n$ for $m=0$) and $\bar Q$
is the algebraic closure of the field of fractions of $R$.  
\end{prop}

\begin{proof}
The inclusion of the
generic point $\Spec \bQ\to \Spec \bZ[1/p]$ induces an isomorphism 
\[
H^{0}_{\et}(\Spec \bZ[1/p],\bZ/p^{n}(k))\to 
H^{0}_{\et}(\Spec \bQ,\bZ/p^{n}(k)),
\]
q.v.~\cite[Proposition~1]{Soule-KtheoryNumRing}.  This reduces to the
case $Q=\bQ$ or $\bQ\phat$ and the \'etale cohomology
$H^{0}_{\et}(\Spec Q;\bZ/p^{n}(k))$ becomes the Galois cohomology
$H^{0}_{\Gal}(Q;\mu_{p^{n}}^{\otimes k}(\bar Q))$. (We will now fix
$\bar Q$ and write $\mu_{p^{n}}$ for $\mu_{p^{n}}(\bar Q)$.) Letting
$G=\Gal(Q(\mu_{p^{n}})/Q)$, the action of $\Gal(\bar Q/Q)$ on
$\mu_{p^{n}}^{\otimes k}$ factors through $G$, and we can identify
$H^{0}_{\Gal}(Q;\mu_{p^{n}}^{\otimes k})$ as the $G$-fixed point
subgroup of $\mu_{p^{n}}^{\otimes k}$.  We have a canonical
isomorphism $G=(\bZ/p^{n})^{\times}$ given by letting $r\in
(\bZ/p^{n})^{\times}$ act on $\alpha \in \mu_{p^{n}}$ by $\alpha\mapsto \alpha^{r}$;
$r$ then acts on $\mu_{p^{n}}^{\otimes k}$ by the $r^{k}$ power map
(i.e., multiplication by $r^{k}$ when we write the group operation
additively). Choosing $r$ to be a generator of $(\bZ/p^{n})^{\times}$,
the $G$-fixed point subgroup of $\mu_{p^{n}}^{\otimes k}$ is the
subset where $r$ acts by the identity, or equivalently, the subset
$\alpha \in \mu_{p^{n}}^{\otimes k}$ such that $\alpha^{r^{k}-1}=1$.  If $p-1$
does not divide $k$, then $r^{k}-1$ is not congruent to $0$ mod $p$,
and the only fixed point is the identity.  On the other hand,
$r^{m(p-1)}-1$ is divisible by $p^{i}$ (and for $i<n$ not $p^{i+1}$)
where $p^{i}=\gcd(|m|p,p^{n})$ (for $m\neq 0$ or $i=n$ if $m=0$), and
the $G$-fixed point subgroup is exactly the subgroup
$\mu_{p^{i}}^{\otimes k}$.
\end{proof}

Defining $H_{\et}^{*}(-;\bZ\phat(k))$ as the inverse limit of
$H_{\et}^{*}(-;\bZ/p^{n}(k))$, we see from the preceding proposition
that for $R=\bZ[1/p]$ or $\bQ\phat$, $H^{0}_{\et}(R;\bZ\phat(k))=0$
for $k\neq 0$.  Since in these cases the homotopy 
groups of $K(R)$ are finitely generated $\bZ\phat$-modules in each
degree (see~\cite[\S 4]{DwyerAX}, \cite[Theorem~D]{HM2}, and~\cite[0.7]{BokstedtMadsen}), 
$\pi_{*}K(R)\phat\iso \lim \pi_{*}(K(R);\bZ/p^{n})$.
Combining these observations and the left exactness of $\lim$, we then
get isomorphisms  
\begin{equation}\label{eq:KvsH}
\begin{aligned}
\pi_{2k-1}(K(R)\phat)&\iso H^{1}_{\et}(\Spec R;\bZ\phat(k))\\
\pi_{2k}(K(R)\phat)&\iso H^{2}_{\et}(\Spec R;\bZ\phat(k+1))
\end{aligned}
\end{equation}
for $k>1$.  Combining these isomorphisms with the fact that
\[
\pi_{2m(p-1)}\hofib(\trc{p}) \iso \pi_{2m(p-1)}\hofib(\cmp) \iso \pi_{2m(p-1)}\hofib(\cmppi)
\]
and $\pi_{2m(p-1)}\hofib(\cmppi)$ fits in an exact sequence
\[
\xymatrix@-1.25pc{%
&&\pi_{2m(p-1)+1}K(\bZ[1/p])\phat\ar[r]
&\pi_{2m(p-1)+1}K(\bQ\phat)\phat
\ar `r/.5pc[d] `d[l] `[llld] `d[r] [lld]
&\\
&\pi_{2m(p-1)}\hofib(\cmppi)\ar[r]
&\pi_{2m(p-1)}K(\bZ[1/p])\phat,
&&
}
\]
Lemma~\ref{lem:main} is now an immediate consequence of the following
pair of lemmas proved in the next section.

\begin{lem}\label{lem:ksurj}
Let $p$ be an odd prime.
The map of rings $\bZ[1/p]\to \bQ\phat$ induces a surjection
\[
H^{1}_{\et}(\Spec\bZ[1/p];\bZ\phat(m(p-1)+1))\to H^{1}_{\et}(\Spec\bQ\phat;\bZ\phat(m(p-1)+1))
\]
for all $k>0$.
\end{lem}

\begin{lem}\label{lem:kvanish}
Let $p$ be an odd prime.
$H^{2}_{\et}(\Spec\bZ[1/p];\bZ\phat(m(p-1)+1))=0$ for all $k>0$.
\end{lem}

We can also deduce Lemma~\ref{lem:nilKZ}:  Quillen's computation of
the $K$-theory of finite fields implies in particular that
$K(\bZ/p)\phat\simeq H\bZ\phat$.  It then follows from Quillen's
localization sequence~\eqref{eq:QuillenLoc} that the map
$K(\bZ)\phat\to K(\bZ[1/p])\phat$ induces an isomorphism in homotopy
groups above degree $1$.  Lemma~\ref{lem:nilKZ} now follows from the
isomorphisms~\eqref{eq:KvsH} and Lemma~\ref{lem:kvanish}.

\section{Proof of Lemmas~\ref{lem:ksurj} and~\ref{lem:kvanish}}

In this section, we prove Lemmas~\ref{lem:ksurj}
and~\ref{lem:kvanish}.  Lemma~\ref{lem:ksurj} is about the
$p$-completion map in \'etale cohomology and the basic tool for
studying this is the Tate-Poitou duality long exact
sequence~\cite{Tate-Duality}.  (Again, for examples applied to
$K$-theory, see~\cite[3.1]{Rognesp}, \cite[\S4]{Geisser-trace},
\cite{Lars-L}.)

In our context, the Tate-Poitou sequence takes the
following form.  Let $M$ be a finite abelian $p$-group with an
action of the Galois group $G$ of the maximal extension of $\bQ$
unramified except at $p$ (e.g., $M=\bZ/p^{n}(k)$) and let $(-)^{*}$
denote the Pontryagin dual, $A^{*}=\Hom(A,\bQ/\bZ)$; then $M^{*}(1)$ is
the $G$-module $\Hom(M,\mu_{\infty})$ where $\mu_{\infty}$ denotes the
$G$-module of all roots of $1$ in $\bar Q$.  The low dimensional part
of Tate-Poitou duality in the case at hand is then summarized by the
following long exact sequence~\cite[3.1]{Tate-Duality}: 
\begin{equation}\label{eq:fintatepoitou}
\begin{gathered}
\xymatrix@-1pc{%
0\ar[r]
&H^{0}_{\et}(\bZ[1/p];M)\ar[r]
&H^{0}_{\et}(\bQ\phat;M)\ar[r]
&(H^{2}_{\et}(\bZ[1/p],M^{*}(1)))^{*}
\ar `r/.5pc[d] `d[l] `[llld] `d/.5pc[lld] [lld]
\\
&H^{1}_{\et}(\bZ[1/p];M)\ar[r]
&H^{1}_{\et}(\bQ\phat;M)\ar[r]
&(H^{1}_{\et}(\bZ[1/p],M^{*}(1)))^{*}
\ar `r/.5pc[d] `d[l] `[llld] `d/.5pc[lld] [lld]
\\
&H^{2}_{\et}(\bZ[1/p];M)\ar[r]
&H^{2}_{\et}(\bQ\phat;M)\ar[r]
&(H^{0}_{\et}(\bZ[1/p],M^{*}(1)))^{*}\ar[r]&0
}
\end{gathered}
\end{equation}
When $M=\bZ/p^{n}(k)$, the first
map in the sequence above 
\[
H^{0}_{\et}(\bZ[1/p];\bZ/p^{n}(k))\to
H^{0}_{\et}(\bQ\phat;\bZ/p^{n}(k))
\]
is an isomorphism by Proposition~\ref{prop:BernoulliDenom}.  Likewise 
when $M=\bZ/p^{n}(k)$ for $k>1$, we see from~\eqref{eq:finKvsH} that
$H^{i}_{\et}(R;\bZ/p^{n}(k))$ is finite for $R=\bZ[1/p]$ or
$\bQ\phat$, 
and it follows that the above is an exact sequence of finite groups.
Taking the inverse limit over $n$ is then exact and we get the
following Tate-Poitou sequence.
\begin{equation}\label{eq:tatepoitou}
\begin{gathered}
\xymatrix@-1pc{%
&&\relax\hphantom{(\bZ[1/p];\bZ\phat(k))}0\ar[r]
&(H^{2}_{\et}(\bZ[1/p],\bZ/p^{\infty}(1-k)))^{*}
\ar `r/.5pc[d] `d[l] `[llld] `d/.5pc[lld] [lld]
\\
&H^{1}_{\et}(\bZ[1/p];\bZ\phat(k))\ar[r]
&H^{1}_{\et}(\bQ\phat;\bZ\phat(k))\ar[r]
&(H^{1}_{\et}(\bZ[1/p],\bZ/p^{\infty}(1-k)))^{*}
\ar `r/.5pc[d] `d[l] `[llld] `d/.5pc[lld] [lld]
\\
&H^{2}_{\et}(\bZ[1/p];\bZ\phat(k))\ar[r]
&H^{2}_{\et}(\bQ\phat;\bZ\phat(k))\ar[r]
&(H^{0}_{\et}(\bZ[1/p],\bZ/p^{\infty}(1-k)))^{*}\ar[r]&0
}
\end{gathered}
\end{equation}

For Lemmas~\ref{lem:ksurj} and~\ref{lem:kvanish}, we
apply~\eqref{eq:tatepoitou} with $k=m(p-1)+1$, combined with the main
theorem of Bayer-Neukirch~\cite{BayerNeukirch}, which relates the
values of the Iwasawa $p$-adic $\zeta$-function with the size of
\'etale cohomology groups.  In the following theorem, $|\cdot|_{p}$
denotes the $p$-adic valuation on $\bQ\phat$, normalized so that
$|p^{n}u|_{p}=p^{-n}$ where $u$ is a unit in $\bZ\phat$.

\begin{thm}[Bayer-Neukirch{~\cite[6.1]{BayerNeukirch}}]\label{thm:BN} 
Let $\zeta_{I}(\omega^{0},s)$ denote the Iwasawa zeta function
of~\cite[5.1]{BayerNeukirch} associated to the trivial character
$\omega^{0}$ and the field $\bQ$.  Let $k=m(p-1)+1$ for $m\neq0$.  If
$\zeta_{I}(\omega^{0},k)\neq 0$ then the groups
$H^{*}_{\et}(\bZ[1/p];\bZ/p^{\infty}(1-k))$ are all finite (zero for
$*\geq 2$) and
\[
|\zeta_{I}(\omega^{0},k)|_{p}=
\frac{\#(H^{0}_{\et}(\bZ[1/p];\bZ/p^{\infty}(1-k)))}%
{\#(H^{1}_{\et}(\bZ[1/p];\bZ/p^{\infty}(1-k)))}
\]
\end{thm}

The following computation of $|\zeta_{I}(\omega^{0},m(p-1)+1)|_{p}$ is
well-known.

\begin{prop}\label{lem:zeta}
For $m\neq0$ and $k=m(p-1)+1$ 
\[
|\zeta_{I}(\omega^{0},k)|_{p}
=|1/(mp)|_{p}.
\]
\end{prop}

\begin{proof}
The Iwasawa zeta function used by
Bayer-Neukirch~\cite[5.1]{BayerNeukirch} depends on a choice of $q\in
\bZ\phat$ with $q\equiv 1 \mod p$.  For the trivial character,
the formula is then
\[
\zeta_{I}(\omega^{0},s)=\frac{p^{\mu_{0}}g_{0}(q^{1-s}-1)}{1-q^{1-s}}
\]
where for $\bQ$ (and any abelian extension of $\bQ$) $\mu_{0}=0$ as a
case of the Iwasawa ``$\mu=0$'' Conjecture proved by
Ferrero-Washington~\cite{FerreroWashington}
(see~\cite[5.3]{BayerNeukirch}) and $g_{0}(x)$ is the characteristic
polynomial of the action of $T\in \bZ\phat[[T]]\iso \Lambda$ on a
$\Lambda$-module denoted as $e_{0}\aM$ in \cite{BayerNeukirch}.  (Here
$\Lambda$ is the Iwasawa Algebra~\cite[7.1]{Washington-1997} for
$\bZ\phat\iso\Gamma<\Gal(\bQ(\mu_{p^{\infty}})/\bQ)$ with topological generator
$\gamma \leftrightarrow 1+T$ acting by $\alpha\mapsto \alpha^{q}$ for $x\in
\mu_{p^{\infty}}$.)  Since for $k=m(p-1)+1$, $m\neq 0$, 
\[
|1-q^{1-k}|_{p}=|1-q^{-m(p-1)}|_{p}=|1-q^{|m|(p-1)}|_{p}=|1/(m(p-1)p)|_{p}=|1/(mp)|_{p},
\]
it suffices to show that $g_{0}(x)=1$.  
This is a special case of the Main Conjecture of Iwasawa
Theory~\cite[\S6~Conjecture]{MazurWiles} for the trivial character.
Though the exposition preceding~\cite[\S9~Theorem]{MazurWiles} makes
the statement appear ambiguous in the case of the trivial character,
this case was known at least as far back as~\cite{Greenberg}, as we
now discuss for the benefit of those (like the authors) not expert in
this theory.

Washington~\cite[15.37]{Washington-1997} denotes $e_{0}\aM$ as
$\epsilon_{0}\aX$ and $\epsilon_{0}\aX_{\infty}$ and shows that
\[
g_{0}(q(1+T)^{-1}-1)=f(T)u(T)
\]
in $\bZ\phat[[T]]$ for $u(T)$ a unit
power series and $f(x)$ the characteristic polynomial of
$\epsilon_{1}X$, where $X$ is the inverse limit of $X_{n}$ and
$X_{n}\iso A_{n}$ is the $p$-Sylow subgroup of the class group of
$\bQ(\mu_{p^{n}})$.  Greenberg~\cite{Greenberg} denotes $X$ as
$X_{K}$, $\epsilon_{1}X$ as $X_{K}^{[1]}$, and defines
$V^{[1]}=\epsilon_{1}X\otimes_{\bZ\phat} \Omega_{p}$ where
$\Omega_{p}=\bar\bQ_{p}$ is the algebraic closure of $\bQ_{p}$.  The
characteristic polynomial of $\epsilon_{1}X$ and $V^{[1]}$ are
therefore equal, and Greenberg~\cite[Corollary~1]{Greenberg} shows
that $V^{[1]}=0$.  Thus, $f(x)=1$ and we conclude that $g_{0}(x)=1$.
(In fact, $\epsilon_{1}X=0$ and $\epsilon_{1}X_{n}=0$ for all $n$ as
can be seen from \cite[6.16,13.22]{Washington-1997} and Nakayama's Lemma.)
\end{proof}

We can also compute $H^{0}_{\et}(\bZ[1/p];\bZ/p^{\infty}(1-k))$ using
Proposition~\ref{prop:BernoulliDenom}, and for $k=m(p-1)+1$ we get
\[
H^{0}_{\et}(\bZ[1/p];\bZ/p^{\infty}(1-k)) =
H^{0}_{\et}(\bZ[1/p];\bZ/p^{\infty}(-m(p-1)))=\mu^{\otimes
(1-k)}_{p^{i}}(\bar \bQ)\iso \bZ\phat/(mp),
\]
where $mp=p^{i}r$ for $r$ relatively prime to $p$, or more concisely,
$p^{i}=|1/(mp)|_{p}$.  The following proposition is now immediate.

\begin{prop}\label{prop:h1zero}
For $m\neq0$ and $k=m(p-1)+1$, $H^{1}_{\et}(\bZ[1/p];\bZ/p^{\infty}(1-k))=0$.
\end{prop}

Combining the previous proposition with the Tate-Poitou
sequence~\eqref{eq:tatepoitou}, the proof of Lemma~\ref{lem:ksurj} is now clear.
For Lemma~\ref{lem:kvanish}, we need the following $K$-theory
computation of Hesselholt-Madsen and B\"okstedt-Madsen. 

\begin{thm}[Hesselholt-Madsen~{\cite[Theorem~D]{HM2}},
B\"okstedt-Madsen~{\cite[0.7]{BokstedtMadsen}}]
For $m>0$, 
\[
\pi_{2m(p-1)}(K(\bQ\phat)\phat)\iso \bZ\phat/(mp).
\]
\end{thm}

\begin{proof}[Proof of Lemma~\ref{lem:kvanish}]
Let $k=m(p-1)+1$.
By the previous theorem and~\eqref{eq:KvsH}, we have
\[
\#(H^{2}_{\et}(\bQ\phat;\bZ\phat(k)))=|1/(mp)|_{p}
\]
and by
Proposition~\ref{prop:BernoulliDenom}, we have
\[
\#((H^{0}_{\et}(\bZ[1/p];\bZ/p^{\infty}(1-k)))^{*})=
\#(H^{0}_{\et}(\bZ[1/p];\bZ/p^{\infty}(1-k)))=|1/(mp)|_{p}.
\]
Because the map 
\[
H^{2}_{\et}(\bQ\phat;\bZ\phat(k))\to 
(H^{0}_{\et}(\bZ[1/p];\bZ/p^{\infty}(1-k)))^{*}
\]
in the Tate-Poitou sequence~\eqref{eq:tatepoitou} is surjective and
the groups are the same finite cardinality, it
must therefore also be injective.   The map 
\[
(H^{1}_{\et}(\bZ[1/p];\bZ/p^{\infty}(1-k)))^{*}\to
H^{2}_{\et}(\bZ[1/p];\bZ\phat(k))
\]
is therefore surjective, and Proposition~\ref{prop:h1zero} then shows
that $H^{2}_{\et}(\bZ[1/p];\bZ\phat(k))$ is zero.
\end{proof}


\bibliographystyle{plain}
\bibliography{thhtc}

\end{document}